\documentclass[graybox,nosecnum]{svmult}
\usepackage[T1]{fontenc}
\usepackage[utf8]{inputenc}
\usepackage[bitstream-charter]{mathdesign}

\usepackage{amsmath,amssymb}
\usepackage{geometry}
\usepackage{hyperref}
\usepackage{soul} 
\usepackage{natbib}
\usepackage[all]{nowidow}
\usepackage{setspace}
\usepackage{titlesec}
\usepackage{hyperref}
\usepackage{csquotes}
\usepackage{pdflscape}
\usepackage{tikz}          
\usetikzlibrary{arrows.meta, decorations.markings}
\usepackage{enumitem}
\hypersetup{
	colorlinks   = true,
	citecolor    = blue,
	urlcolor	 = green!70!blue
}
\usepackage{soul}

\title*{\vspace*{-1cm} When infinity stopped being embarrassing:\\ The doubly infinite series of Pierre Alphonse Laurent and\\ the mathematical rehabilitation of singularities}
\titlerunning{When infinity stopped being embarrassing}
\author{B. Sriraman and N. Karjanto\thanks{\!\!\!\! Corresponding author.}}
\institute{B. Sriraman \at Department of Mathematical Sciences, The University of Montana, Missoula, MT 59812, United States of America \\ \email{sriramanb@mso.umt.edu}
\and N. Karjanto \at Department of Mathematics, College of Natural Science, Sungkyunkwan University, Suwon 16419, South Korea \\ \email{natanael@skku.edu}}
\date{\footnotesize Updated \today}

\begin{document}
\maketitle \abstract{\vspace*{0.25cm}\\ For the better part of a century, isolated singularities, the points at which a complex function fails to be analytic, were treated as pathological obstructions requiring elaborate avoidance strategies. Pierre Alphonse Laurent (1813--1854), a French military engineer who pursued mathematics between port construction projects at Le Havre, ended this avoidance in 1843 by extending Cauchy's Taylor-type theorem to doubly connected (annular) domains, producing the doubly infinite power series that now bears his name. The conceptual shift was profound: negative-power terms in the expansion, far from signaling a breakdown of the formalism, encode precise geometric information about the singularity. Laurent's contribution arrived through an unhappy institutional trajectory---submitted after a prize deadline, subjected to a priority claim by Cauchy that complicated the publication process, and issued in full only posthumously in 1863---yet it became indispensable to every branch of mathematics and mathematical physics that touches on complex function theory. We reconstruct the mathematical problem Laurent solved, place it within Cauchy's analytic program of the 1830s--1840s, examine the institutional failure that prevented publication of the memoir, document the independent parallel proof by Weierstrass (1841, published 1894), and trace the series' subsequent absorption into the standard toolkit via Briot and Bouquet, Weierstrass's systematic treatment, and the residue calculus. Drawing on Laurent's 1843 \textit{Comptes rendus} notice, Cauchy's Academy report of the same year, Joseph Bertrand's memorial notice (read 1860, published 1890), and the peer-reviewed secondary literature (Neuenschwander 1978, 1981; Manning 1975; Bottazzini 1986; Gray 2015), we also analyze the modern philosophical significance of the series, which we term \enquote{exile mathematics}, and survey its reach into perturbation theory, number theory, probability theory, and quantum field theory. Readers familiar with the theorem but not its institutional history will find here a documented account of why a foundational result was withheld from publication for two decades and how it nevertheless achieved canonical status.} \vspace*{0.5cm}


\bigskip
\begin{landscape}
\begin{figure}[h]
\includegraphics[height=\textheight]{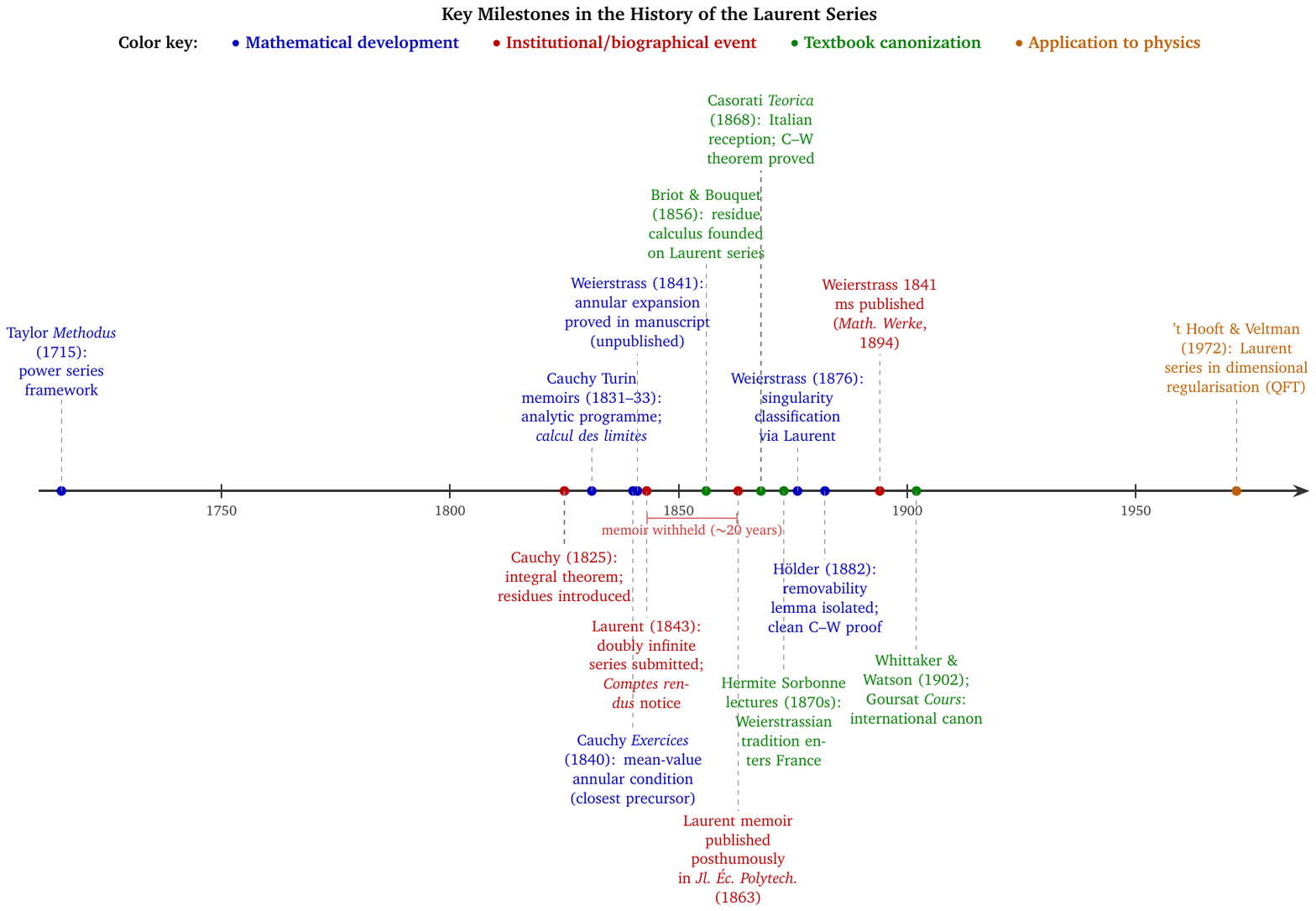}
\end{figure}
\end{landscape}

\section{The mathematical inheritance:\\ Taylor, Cauchy, and the limits of positive powers}			\label{sec:inheritance}

Mathematics has always had a problem with infinity: not the tame infinity of limits approaching a finite value, but the explosive infinity that erupts when a function encounters an isolated singularity. Brook Taylor formalized the standard remedy in the early eighteenth century~\citep{taylor1715methodus}: take any smooth function, examine its derivatives at a base point $c$, and use that local data to build an infinite polynomial approximation. For a function $f$ analytic on the open disk $D(c,R) = \{ z \in \mathbb{C} : |z-c| < R \}$, the Taylor series
\begin{equation}
f(z) = \sum_{n=0}^{\infty} a_n (z - c)^n, \qquad a_n = \frac{f^{(n)}(c)}{n!},		\label{eq:taylor}
\end{equation}
converges absolutely and uniformly on every compact subset of $D(c,R)$. The radius of convergence $R$ is precisely the distance from $c$ to the nearest singularity of $f$~\citep{bottazzini1986higher}. For an entire function like $e^z$ the series converges everywhere:
\[
e^z = \sum_{n=0}^{\infty} \frac{z^n}{n!} = 1 + z + \frac{z^2}{2!} + \frac{z^3}{3!} + \cdots.
\]

The conceptual breakthrough that made a solution possible came from taking seriously the geometry of complex numbers~\citep{needham2023visual}. Throughout the eighteenth century mathematicians had grown increasingly comfortable working with $i = \sqrt{-1}$, initially as a purely algebraic device for solving equations that had no real solutions. What emerged gradually, becoming fully explicit only in the work of Gauss, Argand, and Cauchy in the early nineteenth century, was that complex numbers are not merely numbers but \emph{coordinates in a two-dimensional plane}. Every complex number $z = x + iy$ corresponds to a point whose horizontal displacement from the origin is the real part $x$ and whose vertical displacement is the imaginary part $y$.

This geometric interpretation transformed the study of complex functions. A function $f \colon \mathbb{C} \to \mathbb{C}$ is no longer just a formula to be evaluated: it is a \emph{mapping between planes}, a geometric transformation of the complex plane into itself. Analytic functions, those satisfying the Cauchy--Riemann equations, are precisely the mappings that are locally conformal: they preserve angles and scale uniformly in every direction at non-critical points. From this vantage point, a singularity is not a place where the function \enquote{breaks down}: it is a place where the function does something geometrically distinguished, where the local scaling becomes infinite, where the mapping wraps a neighborhood of the singular point around itself infinitely many times. Singularities, viewed geometrically, are not defects but \emph{features}; they carry precise topological information about the function's global behavior.

Cauchy, working in Turin, then the capital of the Kingdom of Sardinia, where he held a chair at the University after his voluntary exile from France in 1830, produced in his Turin memoirs of 1831--1833 (\enquote{calcul des limites} papers) and the earlier 1825 memoir on definite integrals~\citep{cauchy1825memoire} a geometric theory of complex integration; for Cauchy's rigorous analytical program more broadly, see~\cite{grabiner1981origins}. Cauchy's integral theorem establishes that the integral of a function analytic on and inside a simple closed contour $\gamma$ vanishes identically: $\displaystyle \oint_{\gamma} f(z) \, dz = 0$. When $f$ has isolated singularities inside $\gamma$, this vanishing breaks down in a controlled way described by the \emph{residue theorem}, a consequence of Cauchy's framework that would later find its natural home inside the Laurent expansion (Section~\ref{ssec:residue}):
\begin{equation}
\oint_{\gamma} f(z)\, dz = 2\pi i \sum_{k} \operatorname{Res}(f, \, z_k),			\label{eq:residue-theorem}
\end{equation}
where the sum runs over all singularities $z_k$ enclosed by $\gamma$. Singularities were, in Cauchy's hands, information carriers rather than mere obstacles.

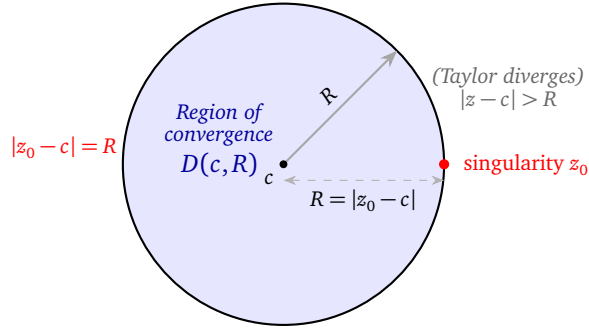
\begin{figure}[h]
\centering
\begin{tikzpicture}[scale=0.85, >=Stealth]
\fill[blue!10] (0,0) circle (2.5cm);
		
\draw[thick, black] (0,0) circle (2.5cm);
		
\draw[->, thick, gray!70](0,0) -- node[above, sloped, black, font=\small]{$R$} (1.77cm, 1.77cm);
		
\node[blue!60!black, font=\small\itshape] at (-1, 0.8){Region of};
\node[blue!60!black, font=\small\itshape] at (-1, 0.4){convergence};
\node[blue!60!black, font=\normalsize] at (-1, 0){$D(c, R)$};
		
\node[gray!80!black, font=\small\itshape] at (3.5, 1){$|z - c| > R$};
\node[gray!80!black, font=\small\itshape] at (3.5, 1.35){(Taylor diverges)};
		
\filldraw[black] (0,0) circle (1.5pt) node[below left, font=\small]{$c$};
		
\filldraw[red] (2.5cm, 0) circle (2pt);
\node[red, font=\small, right=4pt] at (2.5cm, 0) {singularity $z_0$};
		
\node[font=\footnotesize, above left=-1pt] at (-2.5cm, 0){\textcolor{red}{$|z_0 - c| = R$}};
		
\draw[dashed, gray!60, <->] (0,-0.25) -- (2.5cm,-0.25) node[midway, below, font=\footnotesize, black]{$R = |z_0 - c|$};
\end{tikzpicture}
\caption{The disk of convergence $D(c,R) = \{z\in\mathbb{C} : |z-c|<R \}$ (shaded) on which the Taylor series~\eqref{eq:taylor} converges absolutely and uniformly on compact subsets. The radius of convergence $R$ equals the distance from the expansion centre $c$ to the nearest singularity $z_0$ of $f$; inside the disk the series converges; outside it diverges. The structural constraint that the expansion domain must be a \emph{simply connected} disk free of singularities is the limitation that Laurent's theorem removes.}				\label{fig:disk}
\end{figure}

Yet Cauchy's framework contained a structural restriction that he did not fully resolve. Equation~\eqref{eq:taylor} requires $f$ to be analytic on a \emph{simply connected} region containing $c$. As shown in Figure~\ref{fig:disk}, a disk with no holes creates an explicit geometrical constraint. The Taylor series lives inside a disk whose radius is fixed by the distance to the nearest singularity, and it cannot be extended beyond that boundary by any manipulation of the series itself. The moment one attempts to work in a \emph{doubly connected} region---a domain with a hole, such as the annulus $A(c;\,r,R) = \{ z \in \mathbb{C} : r < |z-c| < R \}$, or the punctured neighborhood $0 < |z - c| < \varepsilon$ of an isolated singularity, the Taylor expansion becomes undefined at its foundation. When $c$ is itself a singularity, there is no finite value $f(c)$ from which to form the Taylor coefficients $f^{(n)}(c)/n!$; the very first step of the construction fails. A function as simple as $1/z$ demonstrates the impasse immediately: its only natural expansion point near the origin is the origin itself, yet $1/z$ is not even defined there. The Taylor series is not merely inaccurate in such a region; it is inapplicable.

This limitation was not merely abstract. By the early 1840s, special functions of applied importance (most notably Bessel functions and the elliptic functions of Jacobi) were known to possess isolated singularities, and the problem of representing such functions as series near singular points was pressing~\citep{bottazzini1986higher,gray2015real}. Taylor's theorem offered no traction. Something genuinely new was required.

It is worth delimiting precisely what was new, since negative-power terms in series expansions were not unknown before 1843. Euler, in the \textit{Introductio in analysin infinitorum} (1748) and throughout his operational work, freely manipulated series containing negative powers of the variable, obtaining expansions such as 
\[ 
\frac{1}{1 - x} = -\sum_{n = 1}^{\infty} x^{-n} \qquad \text{for } \; |x| >1  
\]
by algebraic inversion of standard power series. These are genuine negative-power expansions, but they are case-by-case algebraic manipulations: Euler treated power series as what Remmert calls \enquote{non-terminating polynomials}, deploying them as universal formal tools without an integral formula for coefficients, a proof of convergence, or the concept of an annular region as a natural validity domain~\cite[pp.~426--427]{remmert1991theory}. 

The transition from these Eulerian identities to the residue theory of the 1840s required a fundamental shift in perspective: from the expansion of specific algebraic expressions to the representation of any function satisfying the condition of analyticity. What Euler's formal manipulations lacked, and what Cauchy's residue theory extracted only partially, was a systematic framework in which negative-power coefficients could be rigorously determined for an arbitrary analytic function on a prescribed domain. It is precisely this gap that Laurent's 1843 result filled~\citep{grabiner1981origins,remmert1991theory}.

Cauchy, in his theory of residues developed through the 1820s and 1830s, extracted from the singular behavior of a function the single number $\displaystyle \operatorname{Res}(f, \, c) = \frac{1}{2\pi i} \oint_\gamma f(\zeta)\,d\zeta$, which in Laurent's language is the coefficient $a_{-1}$. But Cauchy obtained this one number through his theory of \enquote{extraordinary integrals} without asking what the remaining negative-power coefficients $a_{-2}, a_{-3}, \ldots$ were, or whether they formed part of a systematic series representation valid on a definite domain. In his Turin memoirs of 1831--1833, Cauchy worked with series in annular regions as convergence conditions, but not as the natural habitat of a doubly infinite expansion~\citep{bottazzini1986higher, gray2015real}. The transition from Euler's formal manipulations to Cauchy's regime of precise domains and convergence criteria marks the intellectual context within which Laurent's problem became both formulable and necessary~\cite[pp.~427--428]{grabiner1981origins, remmert1991theory}. What was genuinely new in Laurent's 1843 result was the combination of all three elements that the antecedent literature lacked simultaneously: a \emph{general} representation theorem, valid for any function analytic on \emph{any} annular domain, with \emph{integral formulas for all coefficients} $a_n$, including all negative orders, expressed as contour integrals over any simple closed curve lying in the annulus. It is this combination, and not the isolated appearance of negative powers, that constitutes the advance.

\section{Laurent's discovery: The 1843 memoir and its mathematical content}			\label{sec:discovery}

\subsection{The engineer-mathematician: Pierre Alphonse Laurent}			\label{ssec:bio}

Pierre Alphonse Laurent was born in Paris on 18~July~1813, to an Anglo-French family. He entered the \'{E}cole Polytechnique in~1830, graduated in~1832 as one of the strongest students of his year, and subsequently trained as a military engineer at the \'{E}cole d'Application in Metz, in the Moselle department of northeastern France. After active service in the Algerian campaigns of the 1830s, Algeria then being under French military occupation since the invasion of 1830, including the Tlemcen expedition in northwestern Algeria near the Moroccan border and the Tafna campaign against Abd el-Kader (1837), he returned to France around 1840 and spent approximately six years directing the hydraulic engineering works for the enlargement of the port of Le Havre, in the Seine-Maritime department in the Normandy region of northern France, transforming it into France's principal seaport~\citep{itard1973dsb}.

The biographical fact that matters mathematically is this: Laurent's major result was produced \emph{while he was managing large-scale harbor construction}, not while holding a university post. Joseph Bertrand, in his memorial notice read before the Soci\'{e}t\'{e} des Amis des Sciences in April 1860, is the most nearly contemporaneous scholarly source on Laurent's working conditions: he records that Laurent's superiors held his practical engineering judgment in high regard, which only underlines the intellectual force of the theoretical work he was producing simultaneously while stealing hours from sleep~\cite[pp.~389--393]{bertrand1890}. His professional life and his mathematical research were not two separate activities. The questions engineers asked in the 1840s, about elasticity, wave propagation, heat equilibrium, potential theory, were precisely those whose mathematical formulation involved analytic functions near singularities. Laurent's result did not arrive despite his engineering context; it arrived, in part, because of it.

\subsection{The 1843 memoir and the representation theorem}	\label{ssec:theorem}

In 1843 Laurent submitted a memoir to the Grand Prix competition of the Acad\'{e}mie des Sciences, whose topic had been announced at the Academy's public session of 13~July 1840: to find, in the original French, \textit{les \'{e}quations aux limites que l'on doit joindre aux \'{e}quations ind\'{e}finies pour d\'{e}terminer compl\`{e}tement les maxima et minima des int\'{e}grales multiples} (to find the limiting equations that must be joined to the indefinite equations in order to determine completely the maxima and minima of multiple integrals)~\citep{itard1973dsb, cauchy1844oeuvres}. The prize was awarded to Pierre Fr\'{e}d\'{e}ric Sarrus; a memoir by Delaunay received honorable mention. Cauchy's report on Laurent's submission noted that both Laurent and Sarrus had independently established the same formulas of integration by parts, expressed in different notations, and that Laurent had used methods analogous to those of Sarrus to obtain the limiting equations~\cite[pp.~920--921]{cauchy1844oeuvres}; for a translation of the essential part of this report, see~\citeauthor{todhunter1861} (\citeyear[pp.~476--477]{todhunter1861}).

Laurent's memoir, submitted after the competition closed but before the judges announced their decision, was therefore composed without any knowledge of the competing memoirs~\cite[p.~209]{cauchy1844oeuvres}, a circumstance that underlines the independence of his work. This memoir bears the title \textit{M\'{e}moire sur le calcul des variations}~\cite[p.~62]{itard1973dsb} and is distinct from the separately submitted memoir on the Laurent series theorem, titled \textit{Extension du th\'{e}or\`{e}me de M.~Cauchy relatif \`{a} la convergence du d\'{e}veloppement d'une fonction suivant les puissances ascendantes de la variable}~\cite[p.~938]{cauchy1843rapport}, which was reported to the Academy on 30~October 1843. Both memoirs received Cauchy and Liouville's recommendation for publication in the \textit{Recueil des Savants \'{e}trangers}; neither was published during Laurent's lifetime.

The memoir on the Laurent series theorem was never published; no manuscript copy has been located in the published secondary literature. Its mathematical content survives primarily through Cauchy's report to the Academy of October 1843, which quotes Laurent's own statement of his result. Cauchy's report appears in the \textit{Comptes rendus hebdomadaires} \textbf{17} (1843), pp.~938--940, and it is the closest we have to Laurent's own words. The report is also valuable because it preserves, on the preceding page, Cauchy's own statement of his earlier theorem~\cite[p.~938]{cauchy1843rapport}, which Laurent was explicitly extending. The juxtaposition of the two statements reveals the precise conceptual distance between them. Cauchy's theorem reads, in the original French~\cite[p.~938]{cauchy1843rapport}:
\begin{displayquote}
\textit{$x$ d\'{e}signant une variable r\'{e}elle ou imaginaire, une fonction r\'{e}elle ou imaginaire de $x$ sera d\'{e}veloppable en une s\'{e}rie convergente ordonn\'{e}e suivant les puissances ascendantes de cette variable, tant que le module de la variable conservera une valeur inf\'{e}rieure \`{a} la plus petite de celles pour lesquelles la fonction ou sa d\'{e}riv\'{e}e cesse d'\^{e}tre finie et continue.} 

\medskip

\noindent [Translation: $x$ designating a real or imaginary variable, a real or imaginary function of $x$ will be developable in a convergent series ordered according to the ascending powers of that variable, so long as the modulus of the variable retains a value less than the smallest of those for which the function or its derivative ceases to be finite and continuous.]
\end{displayquote}

Three words carry the entire weight of Laurent's extension: where Cauchy writes \textit{une s\'{e}rie} (one series) ordered according to the \textit{puissances ascendantes} (ascending powers), with the modulus bounded above by a single limit (the geometry of a disk), Laurent writes \textit{la somme de deux s\'{e}ries} (the sum of two series), one ascending and one descending, with the modulus comprised \textit{entre deux limites} (between two limits, the geometry of an annulus). The move from one series to two, from ascending powers only to ascending and descending, and from a disk to an annulus, is the entire theorem. In Laurent's formulation as preserved in Cauchy's
report~\cite[p.~939]{cauchy1843rapport}:
\begin{displayquote}
\textit{$x$ d\'{e}signant une variable r\'{e}elle ou imaginaire, une fonction r\'{e}elle ou imaginaire de $x$ pourra \^{e}tre repr\'{e}sent\'{e}e par la somme de deux s\'{e}ries convergentes, ordonn\'{e}es, l'une suivant les puissances enti\`{e}res et ascendantes, l'autre suivant les puissances enti\`{e}res et descendantes de $x$, tant que le module de $x$ conservera une valeur comprise entre deux limites entre lesquelles la fonction ou sa d\'{e}riv\'{e}e ne cesse pas d'\^{e}tre finie et continue.}

\medskip

\noindent [Translation: $x$ designating a real or imaginary variable, a real or imaginary function of $x$ will be representable by the sum of two convergent series, ordered, one according to the whole and ascending powers, the other according to the whole and descending powers of $x$, so long as the modulus of $x$ retains a value comprised between two limits between which the function or its derivative does not cease to be finite and continuous.]
\end{displayquote}

In the notation that became standard after Weierstrass systematized the theory, what Laurent proved is the following.
\begin{theorem}[Laurent, 1843]			\label{thm:laurent}
Let $f$ be analytic\footnote[1]{The term \textit{holomorphic}, now standard for complex functions that are complex-differentiable on an open domain, was introduced only in 1875 by Briot and Bouquet in the second edition of their 	\textit{Th\'{e}orie des fonctions elliptiques}, as Bottazzini documents~\cite[p.~162]{bottazzini1986higher}. Laurent, Cauchy, and Weierstrass all used the term \textit{analytic} (or, in Cauchy's later terminology, \textit{monogenic} and \textit{synectic}) for what we would now call holomorphic. In the modern statement of the theorem that follows, \textit{analytic} and \textit{holomorphic} are used interchangeably as is standard practice in contemporary mathematics.} on the annulus $A(c;\, r, R) = \{ z \in \mathbb{C} : r < |z - c| < R \}$. Then $f$ admits the unique representation 
\begin{equation}
f(z) = \sum_{n=-\infty}^{\infty} a_n (z - c)^n,			\label{eq:laurent}
\end{equation}
converging absolutely and uniformly on every compact subset of $A(c;\,r,R)$. The coefficients are given by
\begin{equation}
a_n = \frac{1}{2\pi i} \oint_{\gamma} \frac{f(\zeta)}{(\zeta - c)^{n+1}}\, d\zeta,				\label{eq:coeffs}
\end{equation}
where $\gamma$ is any simple closed curve in $A(c;\,r,R)$ encircling $c$ once counterclockwise (winding number $+1$).
\end{theorem}

\begin{remark}
The series~\eqref{eq:laurent} separates into two parts (Figure~\ref{fig:annulus}). The \emph{analytic part} (or regular part) $\displaystyle\sum_{n=0}^{\infty} a_n (z-c)^n$ converges for $|z-c| < R$ and resembles a Taylor series. The \emph{principal part}	$\displaystyle\sum_{n=1}^{\infty} a_{-n} (z-c)^{-n}$ converges for $|z-c| > r$ and encodes the singular behavior. Together they converge on the full annulus. When the inner radius $r = 0$ (a punctured neighborhood of $c$), the	behavior of the principal part determines the nature of the isolated singularity; see Section~\ref{sec:significance}.
\end{remark}
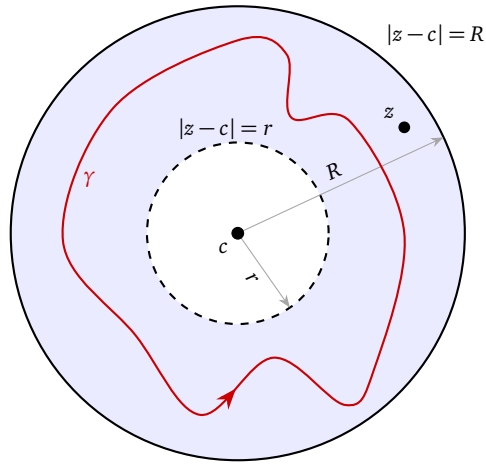
\begin{figure}[h]
\centering
\begin{tikzpicture}[>=Stealth, scale=1]
\fill[blue!8] (0,0) circle (3cm);
\fill[white]  (0,0) circle (1.2cm);

\draw[thick]        (0,0) circle (3cm);
\draw[thick,dashed] (0,0) circle (1.2cm);

\node[above right, font = \small] at ({3*cos(52)}, {3*sin(52)}) {$|z - c| = R$};
\node[below left,  font = \small] at ({1.75*cos(70)}, {1.75*sin(70)}) {$|z - c| = r$};
		
\draw[->, thin, gray!70] (0,0) -- node[above, sloped, font=\small, black]{$R$} ({3*cos(25)}, {3*sin(25)});
\draw[->, thin, gray!70] (0,0) -- node[below, sloped, font=\small, black]{$r$} ({1.2*cos(-55)},{1.2*sin(-55)});
        
\draw[thick, red!80!black, decoration={markings, mark=at position 0.62 with {\arrow[scale=1.2]{Stealth}}}, postaction={decorate}] 
			plot[smooth cycle, tension=0.6] coordinates {
			( 2.20,  0.00)    
			( 1.48,  1.48)    
			( 0.704,  1.5324)  
			( 0.648,  2.3602)  
			( 0.00,  2.55)    
			(-1.65,  1.65)    
			(-2.310,  -0.100)    
			(-1.35, -1.35)    
			(-0.500, -2.40)    
			( 0.470, -1.6446)  
			( 1.421, -2.270)  
			( 1.80, -1.80)    
		};
\node[red!80!black, font=\small] at ({2.082*cos(160)},{2.082*sin(160)}) {$\gamma$};
		
\filldraw[black] (0,0) circle (2.2pt) node[below left, font=\small]{$c$};

\filldraw[black] (2.2,1.4) circle (2pt) node[above left, font=\small]{$z$};
\end{tikzpicture}
\caption{The Laurent annulus $A(c;\,r,R) = \{z \in \mathbb{C} : r < |z-c| < R\}$ (shaded). The solid circle is the outer boundary $|z-c|=R$; the dashed circle is the inner boundary $|z-c|=r$. The contour $\gamma$ (red) is drawn as a general simple closed curve to emphasize that Theorem~\ref{thm:laurent} holds for \emph{any} such curve lying in the annulus and encircling $c$ once counterclockwise with winding number $+1$, not merely for circular contours; a representative point $z$ in the annulus is marked. The coefficients~\eqref{eq:coeffs} are independent of the choice of~$\gamma$.}			\label{fig:annulus}
\end{figure}

To see the theorem at work concretely, consider $f(z) = e^z/z^2$, which has a singularity at the origin. Dividing the known exponential series by $z^2$:
\begin{align*}
\frac{e^z}{z^2}	&= \frac{1}{z^2}\!\left(1 + z + \frac{z^2}{2!}	+ \frac{z^3}{3!} + \cdots\right) \\
				&= \frac{1}{z^2} + \frac{1}{z}	+ \frac{1}{2!} + \frac{z}{3!} + \frac{z^2}{4!} + \cdots.
\end{align*}
The principal part is $z^{-2} + z^{-1}$; the regular part begins at the constant term $\tfrac{1}{2}$. The coefficient $a_{-1} = 1$ is the \emph{residue} of $f$ at $z = 0$.

\subsection{The barrier argument:\\ From doubly connected domains to doubly infinite series}		\label{ssec:proof}

A preliminary remark on the sources is essential. Laurent's 1843 \textit{Comptes rendus} notice~\citep{laurent1843extension} is two pages in length and states the extension theorem without proof. The full memoir on the Laurent series theorem from which that notice was extracted was submitted to the Acad\'{e}mie des Sciences in 1843, reviewed by Cauchy and Liouville, and then, following the institutional failure described in Section~\ref{sec:institution}, never published. It is effectively lost: no manuscript copy has been located in the published secondary literature.\footnote[2]{The Archives de l'Acad\'{e}mie des Sciences (Paris), which holds the \textit{pochettes de s\'{e}ance} and Grand Prix competition papers for 1843, and the Service historique de la D\'{e}fense (Vincennes), which holds the papers of the Corps du G\'{e}nie in which Laurent served, are the repositories where a surviving manuscript would most plausibly be found. A full archival search lies beyond the scope of the present article; the claim of loss follows~\cite[p.~62]{itard1973dsb} and is confirmed by \citeauthor{neuenschwander1981studies2} (\citeyear[p.~88]{neuenschwander1981studies2}) and~\citeauthor{bottazzini1986higher} (\citeyear[p.~178, n.~31]{bottazzini1986higher}), all of whom work exclusively from published sources.} 

The proof technique attributed to Laurent in what follows is therefore a scholarly reconstruction, not a transcription of a documented argument. Its primary basis is Cauchy's report to the Academy of October 1843~\citep{cauchy1843rapport}, which, while written to assert Cauchy's own priority, nonetheless preserves Laurent's statement of his theorem in Laurent's own words and gives some indication of the methods involved. The barrier argument presented below represents the consensus reconstruction in the secondary literature~\citep{neuenschwander1981studies2, bottazzini1986higher, gray2015real}: no published alternative interpretation of Laurent's proof method has been proposed.\footnote[3]{The manuscript presents the barrier argument as the consensus reconstruction of Laurent's proof method. It is worth noting that the history of proofs of Laurent's theorem is itself varied: Remmert documents that Alfred Pringsheim (1896) attempted an integral-free proof via mean values, following the same methodological preference that motivated Cauchy's own 1843 note, while Ludwig Scheeffer (1884) gave a reduction of Laurent's theorem to the Cauchy--Taylor theorem that, though ingenious, Remmert characterises as a \enquote{baroque route} not recommended for teaching since it ultimately depends on the Cauchy integral formula anyway~\cite[pp.~352, 431]{remmert1991theory}. These alternative proof strategies illuminate why the barrier argument, despite its reliance on integration, became the standard pedagogical presentation.} One caveat is necessary: Cauchy's appended note to the same report (p.~940) presents its own derivation of Laurent's theorem via the mean value theorem for analytic functions $f(c) = \frac{1}{2\pi}\int_{-\pi}^{\pi} f(c + re^{ip})\,dp$, which is Cauchy's preferred approach rather than Laurent's. This is documented by~\citeauthor{bottazzini1986higher} (\citeyear[p.~163]{bottazzini1986higher}) and distinguished from Laurent's reconstructed method in~\cite{neuenschwander1981studies2}. The reconstruction that follows reflects Laurent's method as inferred from the structure of his result, not Cauchy's reformulation of it. 

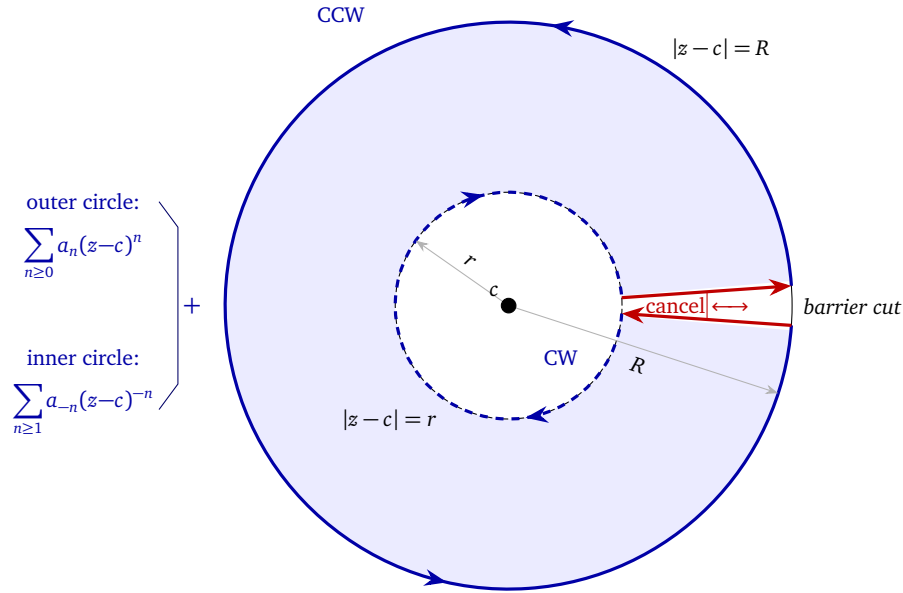
\begin{figure}[h]
\centering
\begin{tikzpicture}[>=Stealth, scale=1.25]	
		
\fill[blue!8] (0,0) circle (3cm);
\fill[white]  (0,0) circle (1.2cm);
\fill[white] (0,0) -- ({3.05*cos(5)},{ 3.05*sin(5)}) arc[start angle=5, delta angle=-10, radius=3.05cm] -- cycle;
		
\draw[thin, black]        (0,0) circle (3cm);
\draw[thin, black, dashed](0,0) circle (1.2cm);
		
\draw[very thick, blue!65!black,
		decoration={markings,
			mark=at position 0.22 with {\arrow{Stealth}},
			mark=at position 0.72 with {\arrow{Stealth}}},
		postaction={decorate}] ({3*cos(4)},{3*sin(4)})
		arc[start angle=4, delta angle=352, radius=3cm];
		
\draw[very thick, blue!65!black, dashed,
		decoration={markings,
			mark=at position 0.22 with {\arrow{Stealth}},
			mark=at position 0.72 with {\arrow{Stealth}}},
		postaction={decorate}] ({1.2*cos(356)},{1.2*sin(356)})
		arc[start angle=356, delta angle=-352, radius=1.2cm];
		
\draw[->, very thick, red!75!black] ({1.2*cos(4)},{1.2*sin(4)}) -- ({3*cos(4)},{3*sin(4)});
\draw[->, very thick, red!75!black]	({3*cos(-4)},{3*sin(-4)}) -- ({1.2*cos(-4)},{1.2*sin(-4)});
\draw[red!75!black, thin] 
		(2.1, { 2.1*tan(4)}) -- (2.1, 0)
		(2.1, {-2.1*tan(4)}) -- (2.1, 0);
\node[red!75!black, font=\footnotesize, right=2pt] at (1.3, 0)	{cancel $\leftarrow\!\!\rightarrow$ };
		
\node[font=\footnotesize, right=1pt] at (3.0, 0) {\textit{barrier cut}};
		
\draw[->, thin, gray!60] (0,0) -- node[below, sloped, font=\small, black]{$R$} ({3*cos(-18)},{3*sin(-18)});
\draw[->, thin, gray!60] (0,0) -- node[above, sloped, font=\small, black]{$r$} ({1.2*cos(145)},{1.2*sin(145)});
		
\node[blue!65!black, font=\small] at ({3.55*cos(120)},{3.55*sin(120)}) {CCW};
\node[blue!65!black, font=\small] at (0.55,-0.55) {CW};
		
\node[font=\footnotesize, align=right, blue!65!black] at (-4.5,  1.1) {outer circle:};
\node[font=\footnotesize, align=right, blue!65!black] at (-4.5,  0.55) {$\displaystyle\sum_{n \ge 0} a_n(z{-}c)^n$};
\node[font=\footnotesize, align=right, blue!65!black] at (-4.5, -0.55) {inner circle:};
\node[font=\footnotesize, align=right, blue!65!black] at (-4.5, -1.1) {$\displaystyle\sum_{n \ge 1} a_{-n}(z{-}c)^{-n}$};
\draw[blue!40!black, thin] (-3.7,  1.1) -- (-3.5,  0.83) -- (-3.5, -0.83) -- (-3.7,-1.1);
\node[blue!65!black, font=\normalsize] at (-3.35, 0) {$+$};
		
\node[font=\small, above right] at ({3*cos(57)},{3*sin(57)})  {$|z - c| = R$};
\node[font=\small, below left]  at ({1.2*cos(237)},{1.2*sin(237)}) {$|z - c| = r$};
		
\filldraw[black] (0,0) circle (2.2pt) node[above left, font=\small]{$c$};
\end{tikzpicture}
\caption{The barrier argument that explains why the Laurent series contains both positive and negative powers. A radial cut (red) connects the inner and outer boundary circles, converting the doubly connected annulus into a simply connected region to which Cauchy's integral formula applies. The contributions along the two sides of the cut traverse opposite directions ($\leftarrow\!\!\rightarrow$) and cancel, since the integrand is analytic there. What remains are two contour integrals: around the outer circle (solid counterclockwise, CCW) and around the inner circle (dashed clockwise, CW). On the outer circle $|\zeta-c|>|z-c|$, so $f(\zeta)/(\zeta-z)$ expands in powers of $(z-c)/(\zeta-c)$, yielding the analytic part (positive powers). On the inner circle $|\zeta-c|<|z-c|$, the expansion runs in powers of $(\zeta-c)/(z-c)$, yielding the principal part (negative powers). Laurent's key insight was that applying this standard Cauchy device to an \emph{annular} domain necessarily produces a doubly infinite series.}	\label{fig:barrier}
\end{figure}

The reconstructed proof runs through a direct application of Cauchy's integral formula to doubly connected domains. The key geometric step is a \emph{barrier argument} (Figure~\ref{fig:barrier}): a radial cut connects the inner and outer boundary circles, converting the doubly connected annulus into a simply connected region. Cauchy's formula is then applied, and the contributions along the two sides of the cut cancel (the integrand is analytic along the cut and the two sides are traversed in opposite directions). This leaves contour integrals along the two circular boundaries. Expanding $f(\zeta)/(\zeta - z)$ geometrically in powers of $(z - c)/(\zeta - c)$ on the outer circle (where $|\zeta-c| > |z - c|$) and in powers of $(\zeta - c)/(z - c)$ on the inner circle (where $|\zeta - c| < |z - c|$) produces the two convergent series whose sum is~\eqref{eq:laurent}, with coefficients~\eqref{eq:coeffs}.

It should be noted that in 1843 neither the language of simply nor doubly connected domains nor the full topological machinery was available: those concepts would be made precise only in Riemann's 1851 dissertation. What Laurent (and Cauchy) worked with was the more informal notion of a region bounded by two curves; the barrier argument was a standard device in Cauchy's own practice. Laurent's innovation was recognizing that applying the barrier argument to an annular region \emph{necessarily} produces a series with both positive and negative powers, not merely the non-negative powers of Taylor's theorem~\citep{neuenschwander1981studies2}. The topological clarification of multiply connected domains and its interaction with the Laurent series in the Weierstrass tradition is traced in detail in~\cite{bottazzini2013hidden}.

\section{The institutional failure} \label{sec:institution}

\subsection{Cauchy's report and priority claim}	\label{ssec:priority}

Cauchy and Liouville~\citep{lutzen1990liouville} were appointed to review Laurent's memoir. Liouville was not a disinterested administrator in this role: Bottazzini documents that he himself began to study doubly periodic functions in the summer of 1844, the year immediately following the review~\cite[p.~163]{bottazzini1986higher}, and his subsequent lectures at the Coll\`{e}ge de France on this subject would later inspire Briot and Bouquet, as documented in Section~\ref{ssec:posthumous}. Both recognized its merit and jointly recommended publication in the \textit{Recueil des savants \'{e}trangers}. What happened next became one of those episodes that punctuate mathematical history. 

The sequence of events is significant: Gray documents that Cauchy made his priority claim at the Academy meeting immediately following Laurent's submission, before he had even submitted the official joint report with Liouville, asserting that earlier results from his \textit{Exercices de math\'{e}matique} of 1826 and his \textit{Exercices d'analyse} of 1840 gave him priority over Laurent~\cite[p.~109]{gray2015real}, \cite[p.~163]{bottazzini1986higher}; that joint report was presented only on 30~October 1843. He appended to the joint report a personal note whose opening sentence reads~\cite[p.~940]{cauchy1843rapport}:\footnote[4]{Some secondary sources, including \citeauthor{bottazzini1986higher}~(\citeyear[p.~163]{bottazzini1986higher}) and \citeauthor{gray2015real}~(\citeyear[p.~109]{gray2015real}), associate 	Cauchy's priority claim with an integral theorem for an annular region from his \textit{Exercices d'analyse et de physique math\'{e}matique} (1840)~\citep{cauchy1840considerations}. A direct examination of that volume (reprinted in \textit{\OE{}uvres compl\`{e}tes}, S\'{e}r.~2, Vol.~11, pp.~331--353) reveals a more nuanced picture. On p.~337--340 Cauchy develops a mean value representation of $f(x)$ as the limit of arithmetic means of $\frac{z}{z-x}f(z)$ at a fixed modulus $r$, valid so long as $f$ and $f'$ remain finite and continuous for moduli between $r_0$ and $R$: this implicitly involves an annular condition, but treats only a single radius at a time and does not produce a doubly infinite series with integral formulas for negative-power coefficients. Remmert confirms this characterization, noting that the 1840 \textit{Exercices} passage is formulated \enquote{without integrals and in terms of mean values}~\cite[pp.~350--351]{remmert1991theory}. The secondary sources are therefore not entirely wrong to point toward this material as relevant prior work, but they overstate its proximity to Laurent's theorem. Cauchy's own appended note to the 1843 report confirms this reading: he states that \enquote{Le th\'{e}or\`{e}me de M.~Laurent peut se d\'{e}duire imm\'{e}diatement} from the 1840 results~\cite[p.~116]{cauchy1843rapport}, \cite[pp.~350--351]{remmert1991theory}, meaning that the simplest route to Laurent's theorem is via the mean value formula rather than via the barrier argument or an explicit annular representation~\cite[p.~163]{bottazzini1986higher}.}
\begin{displayquote}
\textit{Le rapporteur a joint \`{a} ce Rapport la Note suivante, qui indique la mani\`{e}re la plus simple d'arriver au th\'{e}or\`{e}me de M.~Laurent, en partant des principes \'{e}tablis dans les} Exercices d'Analyse et de Physique math\'{e}matique.
\end{displayquote}
In English: the reporter had appended a note indicating the simplest way to arrive at M.~Laurent's theorem starting from principles already established in Cauchy's own earlier \textit{Exercices}, implying that Laurent's result was not original but derivable from existing work. The intellectual authority Cauchy invested in his \textit{Exercices} as a foundation for the whole of analysis is documented by~\cite{grabiner1981origins}.

The claim was, at minimum, difficult to sustain on the mathematical evidence, even though Gray documents that Cauchy had proved a uniqueness theorem for series containing negative powers as early as November 1841~\cite[p.~108]{gray2015real}, \citep{cauchy1841serie}: the existence of such prior work does not bridge the conceptual gap between Cauchy's results, which concerned convergence conditions, and Laurent's theorem, which concerned the representability of an analytic function on an annular domain by a doubly infinite series with explicit integral formulas for all coefficients. It is worth noting that in 1843 Cauchy had not yet arrived at the precise definition of a complex analytic function that his priority claim implicitly required. Gray argues that Cauchy only fully recognized the role of the Cauchy--Riemann equations and the concept of analyticity as late as 1851~\cite[pp.~111--112]{gray2015real}, though Bottazzini's account suggests the recognition was more gradual, developing across the 1840s as Cauchy progressively refined his understanding of what conditions a function must satisfy to admit a power series expansion~\cite[pp.~162--163]{bottazzini1986higher}. On either reading, the conceptual apparatus required to fully justify the priority claim was not in place at the time Cauchy made it in 1843. As~\cite{neuenschwander1981studies2} observes, Laurent's theorem genuinely extended the Cauchyan framework in a non-trivial direction: Cauchy's results concerned series in ascending powers of the variable on \emph{simply connected} domains; Laurent's result concerned doubly infinite series on \emph{annular} domains, which are doubly connected. The formal gap between the two is not a matter of routine reformulation: it requires the barrier argument for doubly connected regions that Cauchy had not made systematic~\citep{gray2015real}. Freudenthal, in his authoritative biography of Cauchy, states plainly that Cauchy had simply failed to discover Laurent's theorem~\cite[p.~141]{freudenthal1971cauchy}, a judgment that aligns with Neuenschwander's and Gray's assessments of the mathematical non-triviality of Laurent's contribution. Kronecker's judgment was more dismissive still: Remmert records that he remarked the Laurent expansion was an \enquote{immediate consequence} of Cauchy's integral theorem and therefore \enquote{useless} to credit to any specific author~\cite[p.~351]{remmert1991theory}---a verdict that, whatever its merits as priority adjudication, underscores how thoroughly the result had been absorbed into the standard apparatus of complex analysis by the late nineteenth century.

Historians have not been unanimous in their assessment of this episode. \cite{neuenschwander1981studies2} and \cite{gray2015real} emphasize the mathematical non-triviality of Laurent's extension and the unjustifiability of Cauchy's priority claim on strictly mathematical grounds. Gray is more direct in his assessment, describing Cauchy's response as characteristic of his habit of taking the work of others as an opportunity to promote his own~\cite[p.~109]{gray2015real}. Belhoste's biography of Cauchy~\citep{belhoste1991cauchy}, however, situates the episode within the norms of the Acad\'{e}mie des Sciences in the 1840s, where senior mathematicians routinely appended priority notes to referees' reports and where such behavior was not regarded as improper by the standards of the institution; the Academy's refereeing and publication procedures in this period are documented in detail by~\cite{grattanguinness1990convolutions}. On this reading, Cauchy was acting within accepted practice, however damaging the consequences for Laurent. \cite{bottazzini1986higher} further notes that the Academy's publication practices in this period were slow and unpredictable regardless of any individual intervention, and that Laurent's late submission, which arrived after the competition deadline, had no formal standing under the Academy's own rules, complicating any simple causal account of why the memoir was never published. Freudenthal's analysis of Cauchy's habitual working method is instructive here: he documents that Cauchy routinely used his role as referee to republish authors' results in broadened form, noting all his own prior related work in the report, behavior he characterizes as standard rather than exceptional for Cauchy~\cite[pp.~133--134]{freudenthal1971cauchy}. The full picture is therefore one of overlapping institutional constraints rather than a single decisive act: Cauchy's priority note, the late submission, the Academy's publication procedures, and the loss of the manuscript all contributed to the outcome. 

Nonetheless, with Cauchy's intervention having complicated the priority question while he simultaneously served as primary referee, the outcome was paradoxical: the Academy formally adopted the conclusions of the joint report~\cite[p.~939]{cauchy1843rapport}, accepting the recommendation to insert Laurent's memoir in the \textit{Recueil des Savants \'{e}trangers}, yet no publication followed. The complete memoir was lost. What Laurent published in his own lifetime amounts to a single two-page notice stating the extension without proof~\citep{laurent1843extension}. The full memoir, with its proofs and elaborations, never appeared in print. Gray notes that when the same misfortune befell a second paper, Laurent turned to applied mathematics, pursuing research in mathematical physics until his death in 1854~\cite[p.~109]{gray2015real}.

The irony is sharp: Cauchy published more than 800 papers~\cite[p.~134]{freudenthal1971cauchy} and became one of the most celebrated mathematicians of his era,\footnote[5]{The figure of 789 papers circulates widely in online sources and popular accounts, but no peer-reviewed scholarly source has been identified as its origin. Freudenthal's authoritative biography in the \textit{Dictionary of Scientific Biography} gives \enquote{more than 800 papers}~\cite[p.~134]{freudenthal1971cauchy}; Belhoste's biography similarly gives \enquote{nearly 800 research articles and treaties}~\cite[p.~vii]{belhoste1991cauchy}. The \textit{\OE{}uvres compl\`{e}tes} edition itself, begun in 1882, remains the primary documentary basis for any count, and its final volume was still outstanding at the time of Freudenthal's writing. Readers should treat any precise figure with appropriate caution.} yet it is Laurent's name, not Cauchy's, that every student of complex analysis encounters first. Bertrand, writing in retrospect, captured the asymmetry: where Cauchy's genius lay in part in his extraordinary productivity, what \enquote{distinguished Cauchy \ldots\ was that he informed the public of them minutely and at length}~\cite[p.~111]{gray2015real}, \cite[pp.~389--393]{bertrand1890}. The Laurent series, the Laurent expansion theorem, Laurent polynomials: these concepts bear his name because his contribution genuinely extended mathematical capability in a way that could not be ignored, whatever the institutional obstacles it encountered.

\subsection{Posthumous publication and canonical reception}	\label{ssec:posthumous}

Laurent died in Paris on 2 September 1854, aged 41,\footnote[6]{\citeauthor{itard1973dsb} (\citeyear[p.~63]{itard1973dsb}) gives the age as forty-two; the arithmetic from the confirmed birth date of 18~July~1813 gives 41~years and 46~days.} leaving a wife and three children~\cite[p.~63]{itard1973dsb}. \citeauthor{bertrand1890} (\citeyear[pp.~392--393]{bertrand1890}) attributes his death to overwork sustained over years of simultaneous professional and mathematical labor. After his death, his widow arranged for two further memoirs to be presented to the Academy. The memoir on light and wave theory was reviewed by Cauchy with a recommendation for publication that the Academy again ignored. A second memoir, \textit{M\'{e}moire sur la th\'{e}orie des imaginaires, sur l'\'{e}quilibre des temp\'{e}ratures et sur l'\'{e}quilibre d'\'{e}lasticit\'{e}}, eventually appeared in the \textit{Journal de l'\'{E}cole Polytechnique} in 1863, nine years after Laurent's death~\citep{laurent1863memoire}. A further document bears on Laurent's posthumous reception: Cauchy submitted a report to the Academy in 1855 reviewing two of Laurent's memoirs, published in the \textit{Comptes rendus} \textbf{40} (1855), pp.~632--634, as noted by~\citeauthor{todhunter1861}~(\citeyear[p.~477]{todhunter1861}).

The theorem that carries his name entered the standard literature not through any work Laurent himself was able to publish in complete form, but through a specific chain of transmission documented by~\cite{neuenschwander1981studies2}: a sequence of textbooks and lecture courses that progressively installed the Laurent series as a foundational tool of French complex analysis and secured its \emph{canonical status} over the four decades following his death.

The chain begins with Briot and Bouquet's 1856 memoir \textit{\'{E}tude des fonctions d'une variable imaginaire}, published in the \textit{Journal de l'\'{E}cole Polytechnique} (Cahier~36, pp.~85--131)~\citep{briotbouquet1856}, in which the residue calculus was for the first time systematically derived from Laurent expansions and the classification of isolated singularities made to depend on the structure of the principal part~\cite[pp.~90--91]{neuenschwander1981studies2}. This was the first instance in the published literature in which Laurent's theorem appeared, not as an isolated result, but as a structural foundation for a broader theory. The same authors returned to the material in their monograph \textit{Th\'{e}orie des fonctions elliptiques} (2nd~ed., 1875)~\citep{briotbouquet1875}, which reached a wider audience, introduced the term \textit{holomorphic} for complex-differentiable functions~\cite[p.~162]{bottazzini1986higher}, and became a standard reference for the next generation of French analysts.

Briot and Bouquet had attended Liouville's lectures on doubly periodic functions at the Coll\`{e}ge de France in 1850--51~\cite[p.~173]{bottazzini1986higher}, \citep{briotbouquet1859}, a circumstance that helps explain the depth of their engagement with the Laurent series as a structural tool. Their book remained the standard text of the French school, was widely adopted in Italy and Germany~\cite[p.~88]{neuenschwander1981studies2}, and was translated into German in 1862~\cite[p.~174]{bottazzini1986higher}, providing a further channel through which the Laurent series entered the German mathematical community alongside Weierstrass's Berlin lectures.

From Briot and Bouquet, the theorem passed into the lecture courses of Charles Hermite at the Sorbonne during the 1870s and 1880s. Hermite’s role in this transmission was more complex than a simple continuation of the French tradition. As Sinkevich documents, Hermite considered himself a student of Weierstrass and explicitly stated that the goal of his Sorbonne lectures was to provide French students with an account of Weierstrass's works and advances~\cite[pp.~85--86]{sinkevich2016weierstrass}. Hermite thus served as the primary conduit between the Weierstrassian algebraic tradition and the French pedagogical pipeline, ensuring that the Laurent series entered the French curriculum with the rigorous power-series foundations that Weierstrass had developed in Berlin, rather than through the Cauchyan integration-theoretic approach alone. These lectures and Hermite's notes, circulated among students in lithographed form by \'{E}mile Picard and subsequently published~\citep{hermite1873cours}, carried this synthesized Laurent series into the training of the generation that would include Poincar\'{e} and Hadamard. The process of canonization was further solidified by Camille Jordan’s \textit{Cours d'analyse de l'\'{E}cole Polytechnique}~\citep{jordan1883cours} (3~vols., 1882--1887), in which the Laurent series appears as a standard instrument of complex function theory, stated, proved, and applied without special remark. Sinkevich further identifies Goursat~\cite[p.~86]{sinkevich2016weierstrass} as being among the mathematicians whose work was shaped by this Weierstrassian influence via Hermite, a lineage that helps explain why Goursat's own 1911 \textit{Cours} could so effectively bridge the two traditions, explicitly citing Cauchy's 1843 result while maintaining a rigorous modern presentation~\cite[p.~92]{goursat1911cours}. In the forty years since his death, Laurent's theorem had thus traveled from an unpublished competition memoir to a textbook commonplace.

The canonization documented above followed primarily the French pedagogical pipeline. Two further channels completed the theorem's absorption into the international mathematical curriculum. In Italy, Casorati's \textit{Teorica delle funzioni di variabili complesse} (1868)~\citep{casorati1868teorica} drew directly on the Briot--Bouquet tradition and introduced the Laurent series to Italian mathematicians; Casorati's work is all the more significant given that the same volume contains his independent proof of what would become the Casorati--Weierstrass theorem (Section~\ref{ssec:classification}), which itself depends on the Laurent classification of isolated singularities. The adoption of Briot and Bouquet's textbook in Italy, noted by~\citeauthor{neuenschwander1981studies2}~(\citeyear[p.~88]{neuenschwander1981studies2}), thus initiated an Italian tradition in complex function theory that ran parallel to the German reception documented in Section~\ref{ssec:traditions}.

In France itself, the next generation of textbooks consolidated Laurent's theorem still further. Appell and Goursat's \textit{Th\'{e}orie des fonctions alg\'{e}briques} (1895)~\citep{appellgoursat1895theorie} and, most influentially, Goursat's \textit{Cours d'analyse math\'{e}matique} (1902--1923)~\citep{goursat1902cours} carried the Laurent series into the twentieth-century university curriculum across France and beyond. Goursat's \textit{Cours} in particular became the dominant French analysis textbook of the early twentieth century, disseminating the residue calculus and singularity classification to generations of students who would never have encountered Laurent's name in any other context. In the 1911 edition, Sections~295--296 treat the Laurent series with full integral coefficient formulas, and a footnote on p.~92 explicitly cites Cauchy's 1843 \textit{Comptes rendus} notice as the origin of the annular representation theorem~\cite[pp.~92--98]{goursat1911cours}, thereby preserving in print the direct historical connection that the manuscript's own Section~\ref{ssec:theorem} documents from primary sources. Goursat's Section~296 further documents Paul Appell's generalization of the Laurent representation to regions bounded by any number of circular arcs~\citep{appellgoursat1895theorie}, extending the annular framework toward the topological generality that would later characterize the several-complex-variables setting surveyed in Section~\ref{sec:applications}~\cite[pp.~94--98]{goursat1911cours}. By the time Whittaker and Watson incorporated the same material into \textit{A Course of Modern Analysis} in 1902~\citep{whittaker1927modern}, the international transmission of the theorem was complete.

\section{The Weierstrass parallel}	\label{sec:weierstrass}

\subsection{The 1841 manuscript}	\label{ssec:weierstrass-ms}

The historical record contains a complication that the priority literature has examined carefully. Karl Weierstrass, working in isolation at M\"{u}nster in Westphalia in northwestern Germany in 1841 while employed as a secondary school teacher, proved essentially the same theorem as Laurent (two years earlier) in a paper titled \textit{Darstellung einer analytischen Function einer complexen Ver\"{a}nderlichen, deren absoluter Betrag zwischen zwei gegebenen Grenzen liegt} (\enquote{Representation of an analytic function of a complex variable whose absolute magnitude lies between two given limits}). Bottazzini confirms the independence of the two results, noting that the same theorem had been obtained but not published by Weierstrass two years before Laurent~\cite[p.~178, n.~31]{bottazzini1986higher}; Neuenschwander confirms that Weierstrass's proof of the Laurent theorem in this paper was independent of Cauchy and predates Laurent's own submission~\cite[p.~93]{neuenschwander1981studies2}. The paper is collected in \textit{Mathematische Werke}, vol.~1 (Berlin: Mayer \& M\"{u}ller, 1894), pp.~51--66~\citep{weierstrass1841darstellung}. 

The intellectual occasion for Weierstrass's 1841 manuscript was his study of Christoph Gudermann's course on elliptic and modular functions, which Weierstrass followed in person at M\"{u}nster in 1839 before writing his diploma paper on the subject~\cite[p.~374]{manning1975emergence}. Gudermann's lectures were published in Crelle's \textit{Journal f\"{u}r die reine und angewandte Mathematik} from 1838 to 1843 under the title \textit{Theorie der Modular-Functionen und der Modular-Integrale}~\cite[p.~377]{manning1975emergence}. Gudermann's subject, that is, doubly periodic meromorphic functions on $\mathbb{C}$, is inseparable from the geometry of annular domains. An elliptic function has isolated poles at the vertices of a period lattice; to represent such a function near any one pole one must work in a punctured neighborhood $0 < |z - c| < r$, or more generally in the annular strip $r_1 < |z - c| < r_2$ lying between two consecutive poles. Taylor's theorem is inapplicable at the outset (the natural expansion center is itself a pole), and in 1841 no systematic tool existed for handling such domains.

Weierstrass's response was the manuscript itself: its explicit restriction of the domain of validity to an annular region directly encodes the geometric situation that Gudermann's functions had forced upon him. Manning establishes that it was this encounter with elliptic functions, and not any program of extending Cauchy's integration theorems, that led Weierstrass to the annular representation problem~\cite[pp.~300--305]{manning1975emergence}. The contrast with Laurent is therefore not merely one of national mathematical style (French integration-theoretic versus German algebraic) but one of \emph{motivating problem}: Laurent was generalizing a theorem about contour integration over doubly connected regions; Weierstrass was seeking to represent specific functions whose very definition placed them in such regions. Brizard's modern pedagogical account of elliptic functions illustrates why this motivation is structurally necessary rather than historically contingent: the Weierstrass $\wp$-function arises as the \emph{natural} representation for physical systems governed by cubic energy polynomials precisely because its poles and periodicity demand the annular framework~\citep{brizard2009primer}, a retrospective confirmation of the mathematical compulsion that drove Weierstrass to the annular problem in 1841.

Weierstrass's paper opens with an explicit statement of its domain of validity: the function is assumed defined and continuous for all $x$ with $M < |x| < B$, and from this assumption he derives the doubly infinite series with integral formulas for coefficients formally identical to~\eqref{eq:laurent} and~\eqref{eq:coeffs}. Weierstrass employed Cauchy's integral theorem independently, unaware at the time of Cauchy's own formulation~\cite[p.~362]{manning1975emergence}, whereas Laurent was explicitly extending Cauchy's framework from within it. \cite{manning1975emergence} documents this paper in detail and establishes that Weierstrass did not attempt to publish it. It emerged only in the 1894 \textit{Werke}, forty years after Laurent's death. Sinkevich corroborates this belated publication pattern, noting that Weierstrass's early breakthroughs remained largely unknown to the broader mathematical community for decades because of his consistent preference for communicating through lectures rather than print~\cite[pp.~79--80]{sinkevich2016weierstrass}, a tendency that, as the present article documents, ultimately ceded nomenclatural priority to Laurent despite Weierstrass's earlier independent discovery.

Weierstrass is also known for the 1876 paper \textit{Zur Theorie der eindeutigen analytischen Functionen}~\citep{weierstrass1876zur}, which developed the classification of isolated singularities and the Weierstrass product theorem. These are distinct contributions; the 1876 paper is not the source of the annular expansion theorem, which belongs to the 1841 manuscript~\citep{manning1975emergence,ullrich1994proof}.

\subsection{Two independent traditions and the dissemination of the theorem}	\label{ssec:traditions}

Manning identifies three historically and mathematically distinct approaches to complex analysis, rooted in the work of Cauchy, Weierstrass, and Riemann respectively, each originating from divergent sources~\cite[pp.~358--359]{manning1975emergence}. The present subsection is concerned with two of these traditions; the Riemannian approach, which is treated separately in the footnote below, arises from different geometric motivations and is not directly implicated in the priority question surrounding the Laurent series. The independence of the Laurent--Weierstrass discoveries is the conclusion of both~\cite{manning1975emergence} and~\cite{ullrich1994proof}.

Weierstrass was working in the German algebraic tradition, building on real and complex power series as primary objects; Laurent was working explicitly within the Cauchyan tradition, extending a specific Cauchy theorem about disk convergence to the annular case. Their starting points, methods, and motivations were different, and the formal coincidence of their results was genuine.\footnote[7]{A third tradition, not treated here, is the Riemannian approach to multiply connected surfaces, in which the topology of annular domains arises naturally from branch point structure; see~\cite{bottazzini2013hidden}, ch.~7--8, for a detailed account. Neuenschwander documents that Riemann explicitly	treated the Laurent series in his introductory lectures on complex function theory at G\"{o}ttingen in 1861~\cite[p.~91]{neuenschwander1981studies2}, confirming that the series entered all three traditions simultaneously.} 

The difference extended to how each tradition conceptualized multiply connected domains before Riemann's topological clarification of 1851. In the French Cauchyan tradition, a doubly connected domain was approached through the barrier argument: a geometric cut converted it into a simply connected region, and the contour integral was the primary  analytical tool. The domain was understood as a region bounded by two curves, and its doubly connected character was handled as a practical obstacle to be removed rather than a structural feature to be exploited. In the Weierstrass tradition, by contrast, the annular domain $M < |x| < B$ was specified algebraically as the set of points where a function is defined and continuous between two bounds on the modulus. No geometric cut was needed and no contour integration was invoked: the doubly connected character of the domain was encoded directly in the two-sided inequality, and the series representation was derived from power-series algebra rather than from the deformation of contours~\citep{manning1975emergence}, a contrast documented in broader context by~\cite{bottazzini2013hidden}.

These two ways of thinking about the same geometric object, one through barriers and integrals and the other through modulus bounds and algebraic manipulation, represent genuinely distinct mathematical cultures, though the dichotomy should not be taken as absolute: Cauchy's writings circulated in Germany and influenced the generation of Weierstrass, while Dirichlet's Paris studies in the 1820s carried French analytic methods back into the German tradition through his teaching of Riemann~\citep{bottazzini2013hidden, neuenschwander1981studies2}. The binary is best understood as a difference of primary orientation and motivating problem rather than a strict partition, and its analytical usefulness should not obscure the genuine intellectual exchange that crossed national and linguistic boundaries throughout the period. 

It should be noted that in 1841 Weierstrass was still working under Cauchy's derivative-based concept of an analytic function; the distinctively Weierstrassian definition of analyticity in terms of power series representation emerged only in 1842, as Manning documents~\cite[p.~360]{manning1975emergence}. The binary contrast between the two traditions is therefore most accurate when applied to the mature programs of both mathematicians rather than to the specific moment of the 1841 manuscript. It is therefore unsurprising that when Riemann provided the topological language in 1851 that could describe both orientations simultaneously, the underlying unity of the two approaches became visible for the first time.

The series carries Laurent's name for the sole reason that Laurent, however minimally, published in 1843, while Weierstrass's proof existed only in manuscript until 1894. This is a case where the sociology of mathematics, by which we mean the institutional mechanisms of publication, refereeing, and priority adjudication, determined a nomenclatural outcome that the mathematical content alone could not have predicted. Yet it was Weierstrass who ensured the Laurent series became a \emph{standard tool} rather than a forgotten priority claim. Weierstrass was consumed by the problem of rigor: he wanted to rebuild all of analysis on foundations so unimpeachable that no result could be doubted. Part of this project, as documented in detail by~\cite{manning1975emergence} and~\cite{bottazzini1986higher}, was the systematic classification of all possible behaviors of analytic functions, and the Laurent series was indispensable to that classification.

Weierstrass showed that every function analytic on an annular domain admits a Laurent expansion there, and that the structure of the principal part (whether it is empty, finite, or infinite) determines exhaustively what kind of singularity the inner boundary represents. It was through Weierstrass's lectures at Berlin, delivered repeatedly from 1856 onward~\cite[p.~93]{neuenschwander1981studies2} and disseminated through notes circulated among his students, that the classification of singularities entered the standard toolkit of complex function theory across Europe. Roy suggests that the Laurent expansion itself may never have been explicitly treated in these lectures, Weierstrass possibly being dissatisfied with a proof that relied on integration rather than the power series methods he regarded as conceptually primary~\cite[p.~13]{rodriguez2013complex}; Manning, however, includes the Laurent series among the results Weierstrass communicated through his lectures and informal discussions, and regards it as integral to the Weierstrassian program~\cite[pp.~358--359, 373]{manning1975emergence}. 

The question of whether and how explicitly the Laurent expansion featured in Weierstrass's Berlin lectures thus remains a point of scholarly discussion; what is clear is that the classification of singularities that depended on it became, through those lectures, the common language of complex function theory across Europe. Sinkevich further characterizes Weierstrass's Berlin lectures as constituting an \enquote{individual} method of teaching that built a powerful school of followers through direct oral transmission rather than through published works~\cite[pp.~80--81]{sinkevich2016weierstrass}, which helps explain both the strength of the Weierstrassian influence on the next generation and the difficulty of precisely documenting which results were explicitly covered in any given lecture course. The Laurent series thus reached universal adoption not through the paper Laurent could not publish, nor through the manuscript Weierstrass chose not to submit, but through the pedagogical authority of Weierstrass's systematic program: a further irony in the history of a theorem already rich in them.

\section{Mathematical significance:\\ Singularities, residues, and the geometry of analysis}	\label{sec:significance}

\subsection{The classification of isolated singularities}	\label{ssec:classification}

Theorem~\ref{thm:laurent} is the foundation on which the complete classification of isolated singularities rests. For a function $f$ with an isolated singularity at $c$, expand $f$ in a Laurent series on the punctured disk $0 < |z - c| < \varepsilon$:
\begin{itemize}[leftmargin=1em]
\item Removable singularity: The principal part $\sum_{n=1}^{\infty} a_{-n}(z-c)^{-n}$ vanishes identically ($a_{-n} = 0$ for all $n \geq 1$). Then $f$ can be assigned a finite value at $c$ making it analytic there (Riemann's removability theorem).
	
\item Pole of order $m$: $a_{-m} \neq 0$ and $a_{-n} = 0$ for all $n > m$. Equivalently, $|f(z)| \to \infty$ as $z \to c$ and $\lim_{z \to c}(z-c)^m f(z) = a_{-m} \neq 0$. The function $e^z/z^2$ of Section~\ref{ssec:theorem} has a pole of	order two at the origin.
	
\item Essential singularity: Infinitely many of the $a_{-n}$ are nonzero. Consider 
	\[
	e^{1/z} = \sum_{n=0}^{\infty} \frac{1}{n!\, z^n}	= 1 + \frac{1}{z} + \frac{1}{2!\,z^2} + \cdots,
	\]
which has an essential singularity at the origin. 
\end{itemize}

By the Casorati--Weierstrass theorem, proved independently by Casorati and by Sokhotski\v{i} in 1868, and possessed by Weierstrass since his 1863 Berlin lectures as confirmed by Herman~Schwarz, a student of Weierstrass, in testimony to Casorati in 1880, with a published proof by Weierstrass in 1876~\citep{neuenschwander1978studies1}---Otto H\"{o}lder subsequently gave a cleaner proof in 1882 in which he isolated the key lemma~\citep{holder1882beweis}---the image of any punctured neighborhood of an essential singularity is dense in $\mathbb{C}$. Picard's theorem~\citep{picard1879sur} sharpens this: the image omits at most one point. H\"{o}lder's 1882 proof is notable for having been the first to isolate and state explicitly the key lemma that a function remaining bounded near an isolated singularity is analytically continuable there, thereby connecting the theorem directly to the removable singularity case in the Laurent classification~\cite[pp.~139--140]{neuenschwander1978studies1}. The relative obscurity of Casorati's and Sokhotski\v{i}'s proofs in Germany explains both why H\"{o}lder returned to the problem independently and why it was Weierstrass's 1876 proof rather than the earlier ones that became the standard reference.

The Laurent series gives an organized view of this taxonomy: what appears from one analytical vantage point as chaotic or pathological behavior is revealed, through the negative-power coefficients, to have precise geometric structure.

\subsection{The residue as Laurent coefficient:\\ Cauchy's approach and Laurent's embedding}	\label{ssec:residue}

The coefficient $a_{-1}$ in~\eqref{eq:laurent} has a special status. Setting $n = -1$ in~\eqref{eq:coeffs}, so that $(\zeta - c)^{n+1} = (\zeta - c)^{0} = 1$:
\begin{equation}
	a_{-1} = \frac{1}{2\pi i} \oint_{\gamma} f(\zeta)\, d\zeta	= \operatorname{Res}(f, \, c).			\label{eq:residue-def}
\end{equation}
The residue theorem~\eqref{eq:residue-theorem} then follows directly from~\eqref{eq:residue-def} applied to each enclosed singularity in turn. In the form standardized in Briot and Bouquet's textbooks of 1856 and 1859 (documented in~\cite[p.~91]{neuenschwander1981studies2}), this theorem became the engine behind the evaluation of a vast range of definite integrals throughout mathematics and theoretical physics.

The distinction between Cauchy's residue approach and Laurent's extension is worth making explicit for the reader of a history of mathematics journal. Cauchy arrived at residues, and at formula~\eqref{eq:residue-def}, by a different route: through his theory of \enquote{extraordinary integrals} around singularities, without systematic use of doubly infinite series. Laurent's innovation was to show that the residue is \emph{one coefficient in a complete series representation}, namely the $n = -1$ coefficient of the Laurent expansion. This embedding of the residue into a broader representational framework is what gave later mathematicians, especially Weierstrass, the conceptual tools to classify singularities exhaustively and to develop the theory of meromorphic functions.

\subsection{The connection to Fourier analysis and the $\mathcal{Z}$-transform}		\label{ssec:fourier}

Under the substitution $z = c + e^{i\theta}$ (restricting to the unit circle centered at $c$), the Laurent series~\eqref{eq:laurent} becomes 
\begin{equation}
	f(c + e^{i\theta}) = \sum_{n=-\infty}^{\infty} a_n\, e^{in\theta},		\label{eq:fourier}
\end{equation}
which is precisely a Fourier series in $\theta$. The Laurent series thus sits at the intersection of complex function theory and harmonic analysis. This structural relationship underlies the $\mathcal{Z}$-transform of digital signal processing, in which a discrete-time sequence $\{x[n]\}_{n=-\infty}^{\infty}$ is encoded as $X(z) = \sum_{n=-\infty}^{\infty} x[n]\, z^{-n}$, a Laurent series in the variable $z^{-1}$.

\subsection{The argument principle}		\label{ssec:argument}

The geometric perspective supplied by the Laurent framework also underlies the argument principle, which relates the topological behavior of a function to its analytic structure. For a meromorphic function $f$ with $N$ zeros and $P$ poles inside a simple closed contour $\gamma$ (each counted with multiplicity): 
\begin{equation}
	\frac{1}{2\pi i} \oint_{\gamma}	\frac{f'(z)}{f(z)}\, dz = N - P.			\label{eq:argument}
\end{equation}
This result, whose proof reduces to reading off residues from the Laurent expansions of $f'/f$ near its poles, transforms the study of analytic functions into a study of their singular points. Functions are not merely formulas to manipulate algebraically: they are geometric objects, maps from one complex plane to another, and the Laurent series provides coordinates for those maps even in their most complicated regions.

\section{Philosophical significance: Exile mathematics}		\label{sec:philosophical}

There is a broader conceptual transformation that the Laurent series exemplifies and that deserves explicit statement in a historical account. For roughly a century after Newton and Leibniz, mathematicians treated singularities as places where their tools failed---edges of the map marked \emph{hic svnt dracones}, here be dragons. Taylor series worked where functions were well-behaved and broke down everywhere else. Mathematicians worked around singularities, using Taylor series in safe regions and waving their hands near singular points.

Laurent's contribution changed this posture entirely. The full expansion~\eqref{eq:laurent} has two parts: the regular part that represents the function's behavior in normal space, the domesticated region where calculus works and derivatives exist; and the principal part that represents the function's behavior in \emph{exile}, that is, in the region around the singularity where normal rules collapse. Together, they form a complete description. Singularities are not defects to be hidden but essential features revealed, through the Laurent expansion, to be \emph{structured} rather than merely pathological: their negative-power coefficients encode precise geometric information that the Taylor framework could not access.

This is what we term \emph{exile mathematics}: the rigorous study of functions at their most unruly, in their most singular regions. The term \emph{exile mathematics} is our own retrospective formulation and is not a category that Laurent, Cauchy, or their contemporaries would have recognized: it is offered as an interpretive lens for making sense of the conceptual shift, not as a claim about how nineteenth-century mathematicians understood their own practice. The Laurent series was the first systematic tool for this enterprise. It would be followed, in the twentieth century, by distribution theory~\citep{lutzen1982prehistory, schwartz1950theorie}, by the study of sheaves and cohomology~\citep{dieudonne1989history}, by renormalization in quantum field theory~\citep{collins1984renormalization}, all of them extensions of the same basic insight that the structure of pathology is itself structured~\citep{gray2008plato, corfield2003philosophy}.

The elegance of the Laurent series lies in its duality. The positive powers build up, approaching the function from below. The negative powers tear down, approaching the function from the singular point outward. Together they meet in the middle, in the annular region $r < |z - c| < R$, and there they give complete information.
\[
f(z) = \underbrace{\cdots +\frac{a_{-2}}{(z-c)^2} + \frac{a_{-1}}{z-c}}_{\text{principal part (exile)}} + \underbrace{a_0 + a_1(z-c) + a_2(z-c)^2 + \cdots}_{\text{regular part}}
\]
The beauty is in the synthesis: the recognition that the regular part and the principal part are not opposites but complements, two aspects of a unified mathematical object that exists in the space between convergence and divergence.

\section{Modern applications and generalizations}	\label{sec:applications}

The reach of the Laurent series is as wide as modern analysis, and its entry into each domain follows the same pattern documented in Sections~\ref{sec:institution} and~\ref{sec:weierstrass}: canonical tools rarely travel through their original papers but through the textbooks and lecture courses of the generation that inherits them.

In \textbf{perturbation theory}, physical problems routinely produce singularities, that is, points where naive series expansions break down. The energy levels of a quantum harmonic oscillator with an anharmonic perturbation can be formally expanded in powers of the coupling $\lambda$:
\begin{equation}
	E_n = \hbar\omega\!\left(n + \tfrac{1}{2}\right) + \lambda E_n^{(1)} + \lambda^2 E_n^{(2)} + \cdots,			\label{eq:perturbation}
\end{equation}
where $\hbar$ is the reduced Planck constant, $\omega$ the oscillator frequency, and $E_n^{(k)}$ the $k$-th order correction to the $n$-th energy level. Although $\lambda$ is treated as small in elementary perturbation theory, the series~\eqref{eq:perturbation} is in fact asymptotic rather than convergent for any fixed $\lambda \neq 0$: its radius of convergence in $\lambda$ vanishes, and the singularities of the energy as a function of $\lambda$ in the complex $\lambda$-plane govern the resummation problem entirely~\citep{bender1999advanced}. Reading off the nature of those singularities (whether poles, branch points, or essential singularities) requires precisely the Laurent analysis that Section~\ref{sec:significance} provides: the structure of the principal part of the energy function near each singularity determines both the type of the singularity and the appropriate resummation scheme.

In \textbf{number theory}, let $\Lambda = \mathbb{Z}\omega_1 + \mathbb{Z}\omega_2$ be a period lattice in $\mathbb{C}$. The Weierstrass $\wp$-function associated to $\Lambda$ has the Laurent expansion
\begin{equation}
	\wp(z) = \frac{1}{z^2} + \sum_{n=1}^{\infty} (2n+1)\, G_{2n+2}(\Lambda)\, z^{2n}, \qquad 0 < |z| < \min_{\omega \in \Lambda\setminus\{0\}}|\omega|, 	\label{eq:weierstrass-p}
\end{equation}
where $G_{2k}(\Lambda) = \displaystyle \sum_{\omega \in \Lambda\setminus\{0\}} \omega^{-2k}$ are the Eisenstein series of the lattice~\citep{silverman2009arithmetic}. The pole at $z = 0$ and the positive-power coefficients together encode deep arithmetic information about elliptic curves and, through the theory of complex multiplication, about class fields.

In \textbf{probability theory}, probability generating functions often have singularities on the unit circle that signal phase transitions~\citep{flajolet2009analytic, stanley1999enumerative}. The Laurent expansion near such singularities, combined with the Fourier-series interpretation~\eqref{eq:fourier}, is the natural tool for extracting asymptotic information from the underlying distributions~\citep{flajolet2009analytic, wilf2006generatingfunctionology}.

In \textbf{quantum field theory}, the renormalization of ultraviolet divergences in perturbative calculations involves Laurent expansions in a dimensional regularisation parameter $\varepsilon$, where spacetime dimension $d = 4 - 2\varepsilon$, a technique introduced independently by~\cite{thooft1972regularization} and by~\cite{bollini1972dimensional} in 1972. The poles at $\varepsilon = 0$ capture the divergences to be subtracted, and the finite remainder is the renormalized amplitude~\citep{collins1984renormalization, peskin1995introduction}. Borcherds makes the Laurent structure explicit from a mathematical standpoint, describing renormalization precisely as the subtraction of the principal part of the Laurent series in $\varepsilon$ to recover a well-defined physical limit~\citep{borcherds2011renormalization}; Brunetti, D\"{u}tsch, and Fredenhagen's algebraic approach to renormalization as the extension of distributions to singular points provides a functional-analytic counterpart to the same insight~\citep{brunetti2009perturbative}. Laurent's negative powers are literally the machinery of renormalization. The reach of Laurent series into modern perturbative physics extends further still: in Quantum Chromodynamics, where perturbative series are asymptotic rather than convergent, Borel transforms and conformal mappings of the Borel plane are used to navigate renormalon singularities whose structure is again determined by the Laurent expansion near the relevant poles~\citep{caprini2018perturbative}, a contemporary instance of the barrier argument logic Laurent himself employed in 1843. For a recent illustration of the same Laurent machinery applied to the renormalization of elastic models in condensed matter physics, see~\cite{chen2024advanced}.

The route by which residue calculus became indispensable to theoretical physics mirrors the mathematical transmission documented in Subsection~\ref{ssec:posthumous}. The critical textbook link was Whittaker and Watson's \textit{A Course of Modern Analysis} (1902; 4th~ed.\ 1927)~\citep{whittaker1927modern}, which presented Laurent expansions, the residue theorem, and the classification of singularities in the form that physicists of the early twentieth century absorbed them. Sommerfeld's lecture courses and their published form in \textit{Partial Differential Equations in Physics}~\citep{sommerfeld1949partial} then carried this apparatus directly into quantum mechanics and electrodynamics. By mid-century, the Laurent series and its residue coefficient $a_{-1}$ were standard instruments of theoretical physics, present in every graduate curriculum, a dissemination path that began with Jordan's \textit{Cours d'analyse} in~1882 and ended in the physics lecture hall eighty years later.

These classical applications all operate within the one-variable complex analytic setting that Laurent inhabited. The twentieth century extended the foundational insight, namely that singular behavior is structured and representable, to richer mathematical settings where the annulus is no longer the natural domain but its structural role is preserved.

In \textbf{several complex variables}, the natural domain of convergence for a doubly infinite power series is no longer an annulus but a \emph{Reinhardt domain}, that is, a domain invariant under independent rotations of each complex coordinate, and the Laurent expansion theorem generalizes to this richer geometric setting, with the topology of convergence domains becoming substantially more intricate than the one-variable case.

In \textbf{algebra and representation theory}, Laurent \emph{polynomials}, that is, finite sums $\sum_{n = -N}^{N} a_n t^n$ in a formal variable $t$, are the building blocks of \emph{loop algebras}, defined as tensor products $\mathfrak{g} \otimes_{\mathbb{C}} \mathbb{C}[t, t^{-1}]$ for a Lie algebra $\mathfrak{g}$. These are the affine analogues of finite-dimensional simple Lie algebras and lie at the heart of the theory of affine Kac--Moody algebras and, through it, of two-dimensional conformal field theory. The negative powers of $t$ in $\mathbb{C}[t, t^{-1}]$ carry exactly the same structural role they carry in Laurent's original series: they encode the behavior on the \enquote{other side} of the singular point.

In \textbf{algebraic geometry}, the formal Laurent series field $k(\!(t)\!)$, the field of fractions of the formal power series ring $k[\![t]\!]$ over a field $k$, is the prototype of a complete discrete valuation field. Laurent series expansions at singular points of algebraic curves are the local data from which global invariants (divisors, differentials, and the Riemann--Roch theorem) are assembled, making Laurent's one-variable construction the local model for all of modern arithmetic geometry. The continuing centrality of local expansions in this setting is illustrated by Greuel and Shustin's survey of plane algebraic curves with prescribed singularities, where Puiseux series provide the canonical local description of singular points whose global deformation properties are under investigation~\citep{greuelshustin2021}.

Each of these generalizations preserves Laurent's foundational insight of 1843: that the systematic description of functions near singularities requires series with both positive and negative powers, and that the coefficients of the negative powers are not anomalies but the most information-rich part of the expansion.

\section{Conclusion}	\label{sec:conclusion}

Pierre Alphonse Laurent solved a concrete and pressing problem while working as a military engineer on the docks of Le Havre: how to represent analytic functions near their singularities. His answer, the doubly infinite series~\eqref{eq:laurent}, was withheld from publication by institutional failure, subject to a priority claim by Cauchy, and published in full only nine years after his death. Yet it proved indispensable.

The Laurent series matters historically because it marks a transition in how mathematicians understood the relationship between regular and singular behavior: from a binary opposition (good regions versus bad) to a unified representational framework in which singularities are informative features. It matters mathematically because it is the foundation of the classification of isolated singularities, the residue theorem, the connection to Fourier analysis, and all subsequent developments in complex function theory. And it matters philosophically because it exemplifies what happens when a mathematician, working in isolation, refuses the intellectual convenience of avoidance and instead asks directly what structure singularities possess.

Hardy observed that the mathematician's patterns, like the painter's or the poet's, must be \emph{beautiful}: that beauty is the first test, and that there is no permanent place in mathematics for ugly work~\citep{hardy1940mathematician}. The Laurent series passes that test: it finds precise structure in the most disordered regions of the complex plane, and makes the infinite tractable by revealing that even at a singularity, the coefficients obey an exact integral law.

Laurent's legacy, as Bertrand noted in 1860, was built not through academic maneuvering or prolific publication but through a single powerful idea that proved necessary. He was, as the history of his series demonstrates, the right person asking the right question at the right moment, even if the moment did not recognize him.

\bigskip
\begin{minipage}{0.99\linewidth}
\centering
An epitaph \\ [0.5ex]
\itshape
Infinity draws near---\\
negative powers reveal\\
\hspace*{1em}what collapse conceals.\\[0.5ex]
Exile and order---\\
singular and smooth unite\\
\hspace*{1em}in one equation.
\end{minipage}

\section*{Primary sources}	\addcontentsline{toc}{section}{Primary Sources}

\begin{enumerate}
\item Bertrand, J. \enquote{Notice sur les travaux du Commandant Laurent, lue en Avril 1860 \`{a} la s\'{e}ance annuelle de la Soci\'{e}t\'{e} des Amis des Sciences.} In \textit{\'{E}loges acad\'{e}miques}. Paris: Gauthier-Villars,	1890, pp.~389--393. [The most nearly contemporaneous biographical source.]
	
\item Cauchy, A.-L. \enquote{Rapport sur un m\'{e}moire de M.~Laurent.} \textit{Comptes rendus} \textbf{17} (1843), pp.~938--940. [The single document that preserves Laurent's own formulation of his theorem; also contains Cauchy's priority note.] Reprinted in \textit{\OE{}uvres compl\`{e}tes d'Augustin Cauchy}, Ser.~1, Vol.~8. Paris: Gauthier-Villars, 1893. 

\item Cauchy, A.-L. \enquote{Rapport sur deux m\'{e}moires de M.~Pierre Alphonse Laurent.} \textit{Comptes rendus hebdomadaires des s\'{e}ances de l'Acad\'{e}mie des Sciences} \textbf{40} (1855), pp.~632--634. [Cauchy's report on Laurent's posthumous memoirs, written after Laurent's death in September 1854; noted in~\cite{todhunter1861}.]

\item Itard, J. \enquote{Laurent, Pierre Alphonse.} In C.~C.~Gillispie (ed.), \textit{Dictionary of Scientific Biography}, Vol.~8. New York: Scribner, 1973, pp.~62--64.

\item Laurent, P.-A. \enquote{Extension du th\'{e}or\`{e}me de M.~Cauchy relatif \`{a} la convergence du d\'{e}veloppement d'une fonction suivant les puissances ascendantes de la variable.} \textit{Comptes rendus hebdomadaires des s\'{e}ances de l'Acad\'{e}mie des Sciences} \textbf{17} (1843), pp.~348--349. [The only mathematical work Laurent published in his lifetime; states the extension theorem without proof.]

\item Laurent, P.-A. \textit{M\'{e}moire sur la th\'{e}orie des imaginaires, sur l'\'{e}quilibre des temp\'{e}ratures et sur l'\'{e}quilibre d'\'{e}lasticit\'{e}}. \textit{Journal de l'\'{E}cole Polytechnique}, Cahier~40, Tome~23 (1863), pp.~75--204. [Posthumous; available at \url{gallica.bnf.fr}.]
	
\item Weierstrass, K. \enquote{Darstellung einer analytischen Function einer complexen Ver\"{a}nderlichen, deren absoluter Betrag zwischen zwei gegebenen Grenzen liegt.} Written 1841; first published in \textit{Mathematische Werke}, Vol.~1. Berlin: Mayer \& M\"{u}ller, 1894, pp.~51--66.
\end{enumerate}


\end{document}